\documentclass[oneside,english]{amsart}
\usepackage[T1]{fontenc}
\usepackage[latin9]{inputenc}
\usepackage{babel}
\usepackage{textcomp}
\usepackage{amstext}
\usepackage{amsthm}
\usepackage{amssymb}
\usepackage[unicode=true,pdfusetitle,
 bookmarks=true,bookmarksnumbered=false,bookmarksopen=false,
 breaklinks=false,pdfborder={0 0 1},backref=false,colorlinks=false]
 {hyperref}

\makeatletter
\numberwithin{equation}{section}
\numberwithin{figure}{section}
 \theoremstyle{definition}
 \newtheorem*{defn*}{\protect\definitionname}
  \theoremstyle{plain}
  \newtheorem{thm}{\protect\theoremname}[section]
  \theoremstyle{plain}
  \newtheorem{prop}{\protect\propositionname}[section]
  \theoremstyle{plain}
  \newtheorem{cor}{\protect\corollaryname}[section]
  \theoremstyle{plain}
  \newtheorem{lem}{\protect\lemmaname}[section]

\makeatother

  \providecommand{\definitionname}{Definition}
  \providecommand{\lemmaname}{Lemma}
  \providecommand{\propositionname}{Proposition}
\providecommand{\corollaryname}{Corollary}
\providecommand{\theoremname}{Theorem}

\begin{document}

\title[An addition formula for the Jacobian theta function]{An addition formula for the Jacobian theta function with applications}

\author{Bing He}

\address{College of Science, Northwest A\&F University \\
Yangling 712100, Shaanxi \\
People's Republic of China}

\email{yuhe001@foxmail.com; yuhelingyun@foxmail.com}

\author{Hongcun Zhai}

\address{Department of Mathematics, Luoyang Normal University\\
 Luoyang 471934, People's Republic of China}

\email{zhai$_{-}$hc@163.com}
\begin{abstract}
Liu {[}Adv. Math. 212(1) (2007), 389\textendash 406{]} established
an addition formula for the Jacobian theta function by using the theory
of elliptic functions. From this addition formula he obtained the
Ramanujan cubic theta function identity, Winquist's identity, a theta
function identities with five parameters, and many other interesting
theta function identities. In this paper we will give an addition
formula for the Jacobian theta function which is equivalent to Liu's
addition formula. Based on this formula we deduce some known theta
function identities as well as new identities. From these identities
we shall establish certain new series expansions for $\eta^{2}(q)$
and $\eta^{6}(q),$ where $\eta(q)$ is the Dedekind eta function,
give new proofs for Jacobi's two, four and eight squares theorems
and confirm several $q$-trigonometric identities conjectured by W.
Gosper. Another conjectured identity on the constant $\Pi_{q}$ is
also settled. The series expansions for $\eta^{6}(q)$ led to new
proofs of Ramanujan\textquoteright s congruence $p(7n+5)\equiv0(\bmod7).$ 

\end{abstract}

\thanks{The first author is the corresponding author. }

\keywords{Addition formula; Elliptic function; Jacobian theta function; Power
of Dedekind's eta function; Two, four and eight squares theorems;
$q$-trigonometric identity}

\subjclass[2000]{11E25, 11F27, 33E05.}
\maketitle

\section{Introduction}

Throughout this paper we take $q=\exp(\pi i\tau),$ where $\text{Im }\tau>0.$
To carry out our study we need the definitions of the Jacobian theta
functions.
\begin{defn*}
Jacobian theta functions $\theta_{j}(z|\tau)$ for $j=1,2,3,4$ are
defined as \cite{L3,WW}
\begin{align*}
 & \theta_{1}(z|\tau)=-iq^{\frac{1}{4}}\sum_{k=-\infty}^{\infty}(-1)^{k}q^{k(k+1)}e^{(2k+1)zi}=2q^{\frac{1}{4}}\sum_{k=0}^{\infty}(-1)^{k}q^{k(k+1)}\sin(2k+1)z,\\
 & \theta_{2}(z|\tau)=q^{\frac{1}{4}}\sum_{k=-\infty}^{\infty}q^{k(k+1)}e^{(2k+1)zi}=2q^{\frac{1}{4}}\sum_{k=0}^{\infty}q^{k(k+1)}\cos(2k+1)z,\\
 & \theta_{3}(z|\tau)=\sum_{k=-\infty}^{\infty}q^{k^{2}}e^{2kzi}=1+2\sum_{k=1}^{\infty}q^{k^{2}}\cos2kz,\\
 & \theta_{4}(z|\tau)=\sum_{k=-\infty}^{\infty}(-1)^{k}q^{k^{2}}e^{2kzi}=1+2\sum_{k=1}^{\infty}(-1)^{k}q^{k^{2}}\cos2kz.
\end{align*}

For brevity, we will use $\vartheta_{j}(\tau),\;\vartheta'_{j}(\tau)$
and $\vartheta''_{j}(\tau)$ to denote $\theta_{j}(0|\tau),\;\theta'_{j}(0|\tau)$
and $\theta''_{j}(0|\tau)$ respectively. For convenience, we will
use the familiar notation
\[
(a;q)_{\infty}=\prod_{n=0}^{\infty}(1-aq^{n}).
\]
With this notation, the Jacobi triple product identity can be written
as \cite{B,GR}
\[
\sum_{n=-\infty}^{\infty}(-1)^{n}q^{n(n-1)/2}z^{n}=(q;q)_{\infty}(z;q)_{\infty}(q/z;q)_{\infty},\quad z\neq0.
\]
The well-known Dedekind eta function $\eta(q)$ is defined by 
\[
\eta(q):=q^{\frac{1}{24}}(q;q)_{\infty}.
\]
Using the Jacobi triple product identity we can deduce the Jacobi
infinite product expressions for theta functions:
\begin{align*}
 & \theta_{1}(z|\tau)=2q^{\frac{1}{4}}\sin z(q^{2};q^{2})_{\infty}(q^{2}e^{2zi};q^{2})_{\infty}(q^{2}e^{-2zi};q^{2})_{\infty},\\
 & \theta_{2}(z|\tau)=2q^{\frac{1}{4}}\cos z(q^{2};q^{2})_{\infty}(-q^{2}e^{2zi};q^{2})_{\infty}(-q^{2}e^{-2zi};q^{2})_{\infty},\\
 & \theta_{3}(z|\tau)=(q^{2};q^{2})_{\infty}(-qe^{2zi};q^{2})_{\infty}(-qe^{-2zi};q^{2})_{\infty},\\
 & \theta_{4}(z|\tau)=(q^{2};q^{2})_{\infty}(qe^{2zi};q^{2})_{\infty}(qe^{-2zi};q^{2})_{\infty}.
\end{align*}
Differentiating the first identity with respect to $z$ and then putting
$z\rightarrow0$ gives 
\begin{equation}
\vartheta'_{1}(\tau)=2q^{\frac{1}{4}}(q^{2};q^{2})_{\infty}^{3}.\label{eq:7-3}
\end{equation}
Setting $z=0$ in the other three results we easily arrive at
\begin{align}
 & \vartheta_{2}(\tau)=2q^{\frac{1}{4}}(q^{2};q^{2})_{\infty}(-q^{2};q^{2})_{\infty}^{2},\label{eq:7-4}\\
 & \vartheta_{3}(\tau)=(q^{2};q^{2})_{\infty}(-q;q^{2})_{\infty}^{2},\nonumber \\
 & \vartheta_{4}(\tau)=(q^{2};q^{2})_{\infty}(q;q^{2})_{\infty}^{2}.\nonumber 
\end{align}

With respect to the (quasi) periods $\pi$ and $\pi\tau,$ we have
\begin{align}
 & \theta_{1}(z+\pi|\tau)=\text{\textminus}\theta_{1}(z|\tau),\quad\theta_{1}(z+\pi\tau|\tau)=-q^{\text{\textminus}1}e^{\text{\textminus}2zi}\theta_{1}(z|\tau),\label{eq:1-1}\\
 & \theta_{2}(z+\pi|\tau)=\text{\textminus}\theta_{2}(z|\tau),\quad\theta_{2}(z+\pi\tau|\tau)=q^{\text{\textminus}1}e^{\text{\textminus}2zi}\theta_{2}(z|\tau),\label{eq:1-2}\\
 & \theta_{3}(z+\pi|\tau)=\theta_{3}(z|\tau),\quad\theta_{3}(z+\pi\tau|\tau)=q^{\text{\textminus}1}e^{\text{\textminus}2zi}\theta_{3}(z|\tau),\label{eq:1-3}\\
 & \theta_{4}(z+\pi|\tau)=\theta_{4}(z|\tau),\quad\theta_{4}(z+\pi\tau|\tau)=-q^{\text{\textminus}1}e^{\text{\textminus}2zi}\theta_{4}(z|\tau).\label{eq:1-14}
\end{align}
We also have the following relations:
\begin{align}
 & \theta_{1}\bigg(z+\frac{\pi}{2}|\tau\bigg)=\theta_{2}(z|\tau),\;\theta_{1}\bigg(z+\frac{\pi\tau}{2}|\tau\bigg)=iB\theta_{4}(z|\tau),\label{eq:1-8}\\
 & \theta_{2}\bigg(z+\frac{\pi}{2}|\tau\bigg)=-\theta_{1}(z|\tau),\;\theta_{2}\bigg(z+\frac{\pi\tau}{2}|\tau\bigg)=B\theta_{3}(z|\tau),\label{eq:1-9}\\
 & \theta_{3}\bigg(z+\frac{\pi}{2}|\tau\bigg)=\theta_{4}(z|\tau),\;\theta_{3}\bigg(z+\frac{\pi\tau}{2}|\tau\bigg)=B\theta_{2}(z|\tau),\label{eq:1-10}\\
 & \theta_{4}\bigg(z+\frac{\pi}{2}|\tau\bigg)=\theta_{3}(z|\tau),\;\theta_{4}\bigg(z+\frac{\pi\tau}{2}|\tau\bigg)=iB\theta_{1}(z|\tau),\label{eq:1-16}
\end{align}
where $B=q^{-1/4}e^{-iz}.$ The following trigonometric expansions
for logarithmic derivatives of theta functions \cite{L3} will be
very useful in the sequel:
\begin{align}
 & \frac{\theta'_{1}}{\theta_{1}}(x|\tau)=\cot x+4\sum_{n=1}^{\infty}\frac{q^{2n}}{1-q^{2n}}\sin(2nx),\label{eq:1-5-1}\\
 & \frac{\theta'_{3}}{\theta_{3}}(x|\tau)=4\sum_{n=1}^{\infty}(-1)^{n}\frac{q^{n}}{1-q^{2n}}\sin(2nx).\label{eq:1-6}
\end{align}

Liu \cite[Theorem 1]{L2} proved an addition formula for the Jacobian
theta function by using the theory of elliptic functions. From this
identity he derived the Ramanujan cubic theta function identity, Winquist's
identity, a theta function identities with five parameters, and many
other interesting theta function identities. See Liu \cite{L5,L1},
Shen \cite{S1,S2} and Whittaker, Watson \cite{WW} for more information
dealing with formulas of theta functions by using elliptic functions.
Liu's addition formula is as follows: 
\end{defn*}
\begin{thm}
\label{t1-1} Let $h_{1}(z)$ and $h_{2}(z)$ are two entire functions
of $z$ which satisfy the functional equations 
\[
h(z+\pi)=-h(z)\quad and\quad h(z+\pi\tau)=-q^{-3}e^{-6iz}h(z).
\]
Then there is a constant $C$ independent of $x$ and $y$ such that
\[
\begin{gathered}(h_{1}(x)-h_{1}(-x))(h_{2}(y)-h_{2}(-y))-(h_{2}(x)-h_{2}(-x))(h_{1}(y)-h_{1}(-y))\\
=C\theta_{1}(x|\tau)\theta_{1}(y|\tau)\theta_{1}(x-y|\tau)\theta_{1}(x+y|\tau).
\end{gathered}
\]
\end{thm}
Actually, Theorem \ref{t1-1} has the following equivalent formula,
which may also be viewed as an addition formula for the Jacobian theta
function $\theta_{1}(z|\tau).$
\begin{thm}
\label{t1-2} Let $h_{1}(z)$ and $h_{2}(z)$ are two entire functions
of $z$ which satisfy the functional equations 
\begin{equation}
h(z+\pi)=h(z)\quad and\quad h(z+\pi\tau)=q^{-2}e^{-4iz}h(z).\label{eq:t-1}
\end{equation}
Then there is a constant $C$ independent of $x$ and $y$ such that
\begin{gather*}
(h_{1}(x)+h_{1}(-x))(h_{2}(y)+h_{2}(-y))-(h_{2}(x)+h_{2}(-x))(h_{1}(y)+h_{1}(-y))\\
=C\theta_{1}(x-y|\tau)\theta_{1}(x+y|\tau).
\end{gather*}
\end{thm}
Theorem \ref{t1-2} includes as special cases many interesting theta
function identities. 
\begin{thm}
\label{t1-6} For all complex numbers $x$ and $y,$ we have
\begin{align}
 & \theta_{2}(2y|2\tau)\theta_{3}(x|\tau/2)-\theta_{2}(2x|2\tau)\theta_{3}(y|\tau/2)=\theta_{1}(x-y|\tau)\theta_{1}(x+y|\tau),\label{eq:1-25}\\
 & \theta_{2}(2y|2\tau)\theta_{4}(x|\tau/2)-\theta_{2}(2x|2\tau)\theta_{4}(y|\tau/2)=\theta_{1}(x-y|\tau)\theta_{1}(x+y|\tau),\label{eq:1-26}\\
 & \theta_{3}(2y|2\tau)\theta_{3}(x|\tau/2)-\theta_{3}(2x|2\tau)\theta_{3}(y|\tau/2)=-\theta_{1}(x-y|\tau)\theta_{1}(x+y|\tau),\label{eq:1-27}\\
 & \theta_{3}(2y|2\tau)\theta_{4}(x|\tau/2)-\theta_{3}(2x|2\tau)\theta_{4}(y|\tau/2)=\theta_{1}(x-y|\tau)\theta_{1}(x+y|\tau).\label{eq:1-28}
\end{align}

\end{thm}
The results in Theorem \ref{t1-6} appear to be new and we have not
found them in the literature. 
\begin{thm}
\label{t2-1} For all complex numbers $x$ and $y,$ we have 
\begin{equation}
\theta_{1}(x|\tau)\theta_{1}(y|\tau)=\theta_{2}(x-y|2\tau)\theta_{3}(x+y|2\tau)-\theta_{2}(x+y|2\tau)\theta_{3}(x-y|2\tau).\label{eq:2-1}
\end{equation}
\end{thm}
The identity (\ref{eq:2-1}) was proved by Lawden \cite[(1.4.7)]{La}
using rearrangements of double infinite series and reproved by Zhai
\cite[(17)]{Z} employing the theory of elliptic functions. 

Making the substitutions $x\rightarrow x+\frac{\pi}{2},\;y\rightarrow y+\frac{\pi}{2}$
in $\eqref{eq:2-1}$ and applying $\eqref{eq:1-2}$, $\eqref{eq:1-3}$
and $\eqref{eq:1-8}$ produces the following result: 
\begin{equation}
\theta_{2}(x|\tau)\theta_{2}(y|\tau)=\theta_{2}(x-y|2\tau)\theta_{3}(x+y|2\tau)+\theta_{2}(x+y|2\tau)\theta_{3}(x-y|2\tau).\label{eq:2-2}
\end{equation}
This identity, which is equivalent to \cite[p.140, (16)]{En}, was
known to Jacobi and is also recorded in \cite[(1.4.9)]{La}. 

From $\eqref{eq:2-1}$ and $\eqref{eq:2-2}$ we can readily deduce
that 
\begin{eqnarray*}
 &  & \theta_{1}(x|\tau)\theta_{1}(y|\tau)+\theta_{2}(x|\tau)\theta_{2}(y|\tau)=2\theta_{2}(x-y|2\tau)\theta_{3}(x+y|2\tau),\\
 &  & \theta_{2}(x|\tau)\theta_{2}(y|\tau)-\theta_{1}(x|\tau)\theta_{1}(y|\tau)=2\theta_{2}(x+y|2\tau)\theta_{3}(x-y|2\tau).
\end{eqnarray*}

Making the substitutions: $x+y\rightarrow2x,\;x-y\rightarrow2y$ in
(\ref{eq:2-1}) we can get

\begin{equation}
\theta_{3}(2x|2\tau)\theta_{2}(2y|2\tau)-\theta_{2}(2x|2\tau)\theta_{3}(2y|2\tau)=\theta_{1}(x-y|\tau)\theta_{1}(x+y|\tau).\label{eq:a-4}
\end{equation}
\begin{thm}
\label{t2-2} For all complex numbers $x,\;y,\;u$ and $v,$ we have
\begin{equation}
\begin{aligned} & \theta_{1}(x-u|\tau)\theta_{1}(x+u|\tau)\theta_{2}(y-v|\tau)\theta_{2}(y+v|\tau)\\
 & \quad-\theta_{1}(y-u|\tau)\theta_{1}(y+u|\tau)\theta_{2}(x-v|\tau)\theta_{2}(x+v|\tau)\\
 & =\theta_{2}(u-v|\tau)\theta_{2}(u+v|\tau)\theta_{1}(x-y|\tau)\theta_{1}(x+y|\tau).
\end{aligned}
\label{eq:2.1}
\end{equation}
\end{thm}
Actually, (\ref{eq:2.1}) is equivalent to \cite[Theorem 1.3]{L4}
via the substitutions $x\rightarrow x+\dfrac{\pi+\pi\tau}{2}$ and
$y\rightarrow y+\dfrac{\pi+\pi\tau}{2}$. 
\begin{thm}
\label{t2-3} For all complex numbers $x,\;y$ and $u,$ we have 
\begin{gather}
\begin{gathered}\theta_{1}(x-u|\tau)\theta_{1}(x+u|\tau)\theta_{2}(2y|2\tau)-\theta_{1}(y-u|\tau)\theta_{1}(y+u|\tau)\theta_{2}(2x|2\tau)\\
=\theta_{2}(2u|2\tau)\theta_{1}(x-y|\tau)\theta_{1}(x+y|\tau),
\end{gathered}
\label{eq:22-4}\\
\begin{gathered}\theta_{1}(x-u|\tau)\theta_{1}(x+u|\tau)\theta_{3}(2y|2\tau)-\theta_{1}(y-u|\tau)\theta_{1}(y+u|\tau)\theta_{3}(2x|2\tau)\\
=\theta_{3}(2u|2\tau)\theta_{1}(x-y|\tau)\theta_{1}(x+y|\tau).
\end{gathered}
\label{eq:22-8}
\end{gather}
\end{thm}
\begin{thm}
\label{t2-4} For all complex numbers $x$ and $y$ we have
\begin{equation}
\begin{gathered}\theta_{1}^{3}(x|\tau)\theta_{1}(y|\tau/3)-\theta_{1}^{3}(y|\tau)\theta_{1}(x|\tau/3)\\
=\frac{\vartheta'_{1}(\tau/3)}{\vartheta'_{1}(\tau)}\theta_{1}(x|\tau)\theta_{1}(y|\tau)\theta_{1}(x-y|\tau)\theta_{1}(x+y|\tau).
\end{gathered}
\label{eq:1-20}
\end{equation}
\end{thm}
To the best of our knowledge, these identities (\ref{eq:22-4}), (\ref{eq:22-8})
and (\ref{eq:1-20}) seem to be new. 

The rest of this paper is organizied as follows. We will show the
equivalence between Theorems \ref{t1-1} and $\ref{t1-2}$ in Section
\ref{sec:Equi}. Section \ref{sec:Proof3} is devoted to proof of
Theorem \ref{t1-6} and derivation of several new series representations
of $\eta^{k}(q)$ for $k=2,\:6$. A new proof of Ramanujan's partition
congruence $p(7n+5)\equiv0\pmod7$ is also given. In Section \ref{sec:New},
we prove Theorem \ref{t2-1} using Theorem $\ref{t1-2}.$ At the same
time, we derive another series representation for $\eta^{6}(q)$ using
(\ref{eq:a-4}). In Section \ref{sec:sn}, we shall give new proofs
for Jacobi's two, four and eight squares theorems. In Section \ref{sec:Proof1}
we will prove Theorems \ref{t2-2}\textendash \ref{t2-4} using Theorem
$\ref{t1-2}$. In addition, we employ Theorems $\ref{t2-2}$ and $\ref{t2-3}$
to confirm certain $q$-trigonometric identities conjectured by W.
Gosper. In fact, we will show two general and interesting $q$-trigonometric
identities\footnote{See Section \ref{sec:Proof1} for the definitions of $q$-trigonometric
functions.}:
\[
\begin{aligned} & \sin{}_{q}(x-u)\sin{}_{q}(x+u)\cos_{q}(y-v)\cos_{q}(y+v)\\
 & \quad-\cos_{q}(x-v)\cos_{q}(x+v)\sin{}_{q}(y-u)\sin{}_{q}(y+u)\\
 & =\cos_{q}(u-v)\cos_{q}(u+v)\sin{}_{q}(x+y)\sin{}_{q}(x-y)
\end{aligned}
\]
and
\[
\cos''_{q}0\cos{}_{\sqrt{q}}(2z)=(\sin'_{q}z)^{2}-\sin{}_{q}z\sin_{q}''z+2\cos''_{\sqrt{q}}0\cos_{q}^{2}z.
\]
Another conjectured identity on the constant $\Pi_{q}$ is also proved
by using Theorem \ref{t2-4}.

\section{\label{sec:Equi} The equivalence between Theorems \ref{t1-1} and
$\ref{t1-2}$}

We first prove Theorem $\ref{t1-2}$ using Theorem $\ref{t1-1}.$
Suppose $h_{1}(z)$ and $h_{2}(z)$ satisfy the functional equations
(\ref{eq:t-1}). Take $f_{1}(z)=h_{1}(z)\theta_{1}(z|\tau)$ and $f_{2}(z)=h_{2}(z)\theta_{1}(z|\tau).$
Then $f_{1}(z)$ and $f_{2}(z)$ satisfy the conditions of Theorem
$\ref{t1-1}.$ So
\begin{gather*}
(f_{1}(x)-f_{1}(-x))(f_{2}(y)-f_{2}(-y))-(f_{2}(x)-f_{2}(-x))(f_{1}(y)-f_{1}(-y))\\
=C\theta_{1}(x|\tau)\theta_{1}(y|\tau)\theta_{1}(x-y|\tau)\theta_{1}(x+y|\tau).
\end{gather*}
Applying the identity $\theta_{1}(-z|\tau)=-\theta_{1}(z|\tau)$ and
cancelling out the common factor $\theta_{1}(x|\tau)\theta_{1}(y|\tau),$
we can readily obtain the result of Theorem $\ref{t1-2}.$

We now prove Theorem $\ref{t1-1}$ using Theorem $\ref{t1-2}.$ Suppose
$h_{1}(z)$ and $h_{2}(z)$ satisfy the functional equations of Theorem
$\ref{t1-1}.$ Then $\dfrac{h_{1}(z)-h_{1}(-z)}{\theta_{1}(z|\tau)}$
and $\dfrac{h_{2}(z)-h_{2}(-z)}{\theta_{1}(z|\tau)}$ are two entire
functions and satisfy the functional equations $\eqref{eq:t-1}$.
Set $f_{1}(z)=\dfrac{h_{1}(z)-h_{1}(-z)}{\theta_{1}(z|\tau)}$ and
$f_{2}(z)=\dfrac{h_{2}(z)-h_{2}(-z)}{\theta_{1}(z|\tau)}$ in Theorem
$\ref{t1-2}.$ Then 
\begin{gather*}
(f_{1}(x)+f_{1}(-x))(f_{2}(y)+f_{2}(-y))-(f_{2}(x)+f_{2}(-x))(f_{1}(y)+f_{1}(-y))\\
=4C\theta_{1}(x-y|\tau)\theta_{1}(x+y|\tau).
\end{gather*}
Using the identity $\theta_{1}(-z|\tau)=-\theta_{1}(z|\tau)$ again
and multiplying both sides of the resulting identity by $\theta_{1}(x|\tau)\theta_{1}(y|\tau)/4$
yields the result of Theorem $\ref{t1-1}$ immediately.

\section{\label{sec:Proof3} Proof of Theorem \ref{t1-6} and series representations
of $\eta^{k}(q)$ for $k=2,\:6$ }

We first prove Theorem \ref{t1-6} using Theorem \ref{t1-2}.

\noindent{\it Proof of Theorem \ref{t1-6}.} We first prove (\ref{eq:1-26}).
It follows from (\ref{eq:1-2}) and (\ref{eq:1-14}) that $\theta_{2}(2z|2\tau)$
and $\theta_{4}(z|\tau/2)$ satisfy the conditions of Theorem \ref{t1-2}.
Taking $h_{1}(z)=\theta_{2}(2z|2\tau)$ and $h_{2}(z)=\theta_{4}(z|\tau/2)$
in Theorem \ref{t1-2} we have 
\begin{equation}
4(\theta_{2}(2x|2\tau)\theta_{4}(y|\tau/2)-\theta_{2}(2y|2\tau)\theta_{4}(x|\tau/2))=C\theta_{1}(x-y|\tau)\theta_{1}(x+y|\tau).\label{eq:q}
\end{equation}

We now determine the constant $C$. Applying the infinite product
representations of theta functions we can find that
\[
2\theta_{4}(x-\pi\tau/4|\tau)\theta_{4}(x+\pi\tau/4|\tau)=q^{-1/8}\vartheta_{2}(\tau/2)\theta_{4}(x|\tau/2).
\]
Replacing $x$ by $x+\pi\tau/2$ in this equation and employing (\ref{eq:1-14})
and (\ref{eq:1-16}), we arrive at 
\begin{equation}
2\theta_{1}(x-\pi\tau/4|\tau)\theta_{1}(x+\pi\tau/4|\tau)=q^{-1/8}\vartheta_{2}(\tau/2)\theta_{4}(x|\tau/2).\label{eq:7-5}
\end{equation}
It follows from the definition of $\theta_{2}$ and the Jacobi triple
product identity that 
\begin{equation}
\begin{aligned}\theta_{2}(\pi\tau/2|2\tau) & =q^{\frac{1}{2}}\sum_{k=-\infty}^{\infty}q^{2k(k+1)}e^{(2k+1)\pi i\tau/2}\\
 & =q\sum_{k=-\infty}^{\infty}q^{2k^{2}+3k}\\
 & =q(q^{4};q^{4})_{\infty}(-q^{5};q^{4})_{\infty}(-q^{-1};q^{4})_{\infty}\\
 & =(q^{4};q^{4})_{\infty}(-q;q^{2})_{\infty}.
\end{aligned}
\label{eq:7-6}
\end{equation}
Setting $y=\pi\tau/4$ in (\ref{eq:q}) and using (\ref{eq:7-5})
and (\ref{eq:7-6}) we get $C=-4.$ We substitute $C=-4$ back into
(\ref{eq:q}) to establish the result (\ref{eq:1-26}). 

Making the substitutions: $(1)\;x\rightarrow x+\dfrac{\pi}{2},\;y\rightarrow y+\dfrac{\pi}{2};\;(2)\;x\rightarrow x+\dfrac{\pi+\pi\tau}{2},\;y\rightarrow y+\dfrac{\pi+\pi\tau}{2};\;(3)\;x\rightarrow x+\dfrac{\pi\tau}{2},\;y\rightarrow y+\dfrac{\pi\tau}{2}$
in (\ref{eq:1-26}) we can easily obtain (\ref{eq:1-25}), (\ref{eq:1-27})
and (\ref{eq:1-28}) respectively. This finishes the proof of Theorem
\ref{t1-6}. \qed

Many series representations for $(q;q)_{\infty}^{2}$ were found by
different people in the literature. L.J. Rogers \cite{R} was the
first to prove the following identity
\[
(q;q)_{\infty}^{2}=\sum_{\substack{m,n=-\infty\\
n\geq2|m|
}
}^{\infty}(-1)^{m+n}q^{\frac{n(n+1)}{2}-\frac{m(3m-1)}{2}}.
\]
Ewell in \cite{Ew} found the series expansion for $(q;q)_{\infty}^{2}$
\[
(q;q)_{\infty}^{2}=\sum_{m,n=-\infty}^{\infty}(q^{3m^{2}+3n^{2}+n}-q^{3m^{2}+3n^{2}+3m+2n+1}).
\]
In \cite{S3} Shen deduced the following formula for $(q;q)_{\infty}^{2}$
\[
(q;q)_{\infty}^{2}=\sum_{m,n=-\infty}^{\infty}(-1)^{m}q^{m^{2}+n^{2}+mn+n/2}.
\]
Liu \cite{L4,L7} established two identities involving $(q;q)_{\infty}^{2}$:
\begin{align*}
(q;q)_{\infty}^{2} & =\sum_{n=0}^{\infty}\sum_{-n\leq j\leq n}(-1)^{j}(1-q^{2n+1})q^{2n^{2}+n-j(3j+1)/2},\\
(q;q)_{\infty}^{2} & =\dfrac{1}{2}\sum_{m,n=-\infty}^{\infty}((-1)^{n}-(-1)^{m})q^{(3m^{2}+3n^{2}+4m+1)/4}.
\end{align*}

From Theorem \ref{t1-6} we can derive several series representations
of $\eta^{2}(q).$ 
\begin{prop}
\label{c1} For $|q|<1,$ we have 
\begin{align*}
(q;q)_{\infty}^{2} & =\sum_{m,n=-\infty}^{\infty}(-1)^{m}q^{(3m^{2}+1)/4+3n^{2}}(q^{2n}-q^{m}),\\
(q;q)_{\infty}^{2} & =q\sum_{m,n=-\infty}^{\infty}q^{3m(m+1)+3n^{2}/4}(q^{2m+1}-q^{n}),\\
(q;q)_{\infty}^{2} & =q\sum_{m,n=-\infty}^{\infty}(-1)^{n}q^{3m(m+1)+3n^{2}/4}(q^{2m+1}-q^{n}),\\
(q;q)_{\infty}^{2} & =\sum_{m,n=-\infty}^{\infty}q^{(3m^{2}+1)/4+3n^{2}}(q^{m}-q^{2n}).
\end{align*}
\end{prop}

\noindent{\it Proof.} We first prove the first identity. Applying
the Jacobi infinite product expressions for $\theta_{1}$, we find
that 
\begin{equation}
\theta_{1}(\pi\tau|3\tau)=iq^{-1/4}(q^{2};q^{2})_{\infty}.\label{eq:7-7}
\end{equation}
It follows from the definitions of theta functions that
\begin{equation}
\begin{aligned}\vartheta_{3}(6\tau)=\sum_{k=-\infty}^{\infty}q^{6k^{2}}, & \quad\theta_{3}(2\pi\tau|6\tau)=\sum_{k=-\infty}^{\infty}q^{6k^{2}+4k},\\
\vartheta_{4}(3\tau/2)=\sum_{k=-\infty}^{\infty}(-1)^{k}q^{3k^{2}/2}, & \quad\theta_{4}(\pi\tau|3\tau/2)=\sum_{k=-\infty}^{\infty}(-1)^{k}q^{3k^{2}/2+2k}.
\end{aligned}
\label{eq:7-8}
\end{equation}
Replacing $\tau$ by $3\tau$ in (\ref{eq:1-28}) and then setting
$x=\pi\tau,\;y=0$ we have
\[
\vartheta_{3}(6\tau)\theta_{4}(\pi\tau|3\tau/2)-\vartheta_{4}(3\tau/2)\theta_{3}(2\pi\tau|6\tau)=\theta_{1}^{2}(\pi\tau|3\tau).
\]
Substituting (\ref{eq:7-7}) and (\ref{eq:7-8}) into this equation
and replacing $q^{2}$ by $q$ gives
\[
(q;q)_{\infty}^{2}=\sum_{m,n=-\infty}^{\infty}(-1)^{m}q^{(3m^{2}+1)/4+3n^{2}}(q^{2n}-q^{m}).
\]
This proves the first identity. 

Similarly, from (\ref{eq:1-25})\textendash (\ref{eq:1-27}) we can
also derive the second, third and forth identities. This completes
the proof of Proposition \ref{c1}. \qed

Certain series representations for $(q;q)_{\infty}^{6}$ appeared
in the literature. Schoeneberg \cite{Sc} gave an interesting formula:
\[
(q;q)_{\infty}^{6}=\sum_{m,n=-\infty}^{\infty}\Re(m+2ni)^{2}q^{(m^{2}+4n^{2}-1)/4},
\]
where $\Re(z)$ denotes the real part of the complex number $z.$ 

In \cite{Ew} Ewell also established a series expansion for $(q;q)_{\infty}^{6}.$
Liu in \cite{L4} found the series expansion for $(q;q)_{\infty}^{6}:$
\[
(q;q)_{\infty}^{6}=\dfrac{1}{4}\sum_{m,n=-\infty}^{\infty}(-1)^{m}(n^{2}-m^{2})q^{(m^{2}+n^{2}-1)/4}.
\]

We now give several new series representations for $\eta^{6}(q)$
using Theorem \ref{t1-6}.
\begin{prop}
\label{c2} For $|q|<1,$ we have
\begin{align}
(q;q)_{\infty}^{6} & =\dfrac{1}{2}\sum_{m,n=-\infty}^{\infty}((2m+1)^{2}-n^{2})q^{m(m+1)+n^{2}/4},\nonumber \\
(q;q)_{\infty}^{6} & =\dfrac{1}{2}\sum_{m,n=-\infty}^{\infty}(-1)^{n}((2m+1)^{2}-n^{2})q^{m(m+1)+n^{2}/4},\nonumber \\
(q;q)_{\infty}^{6} & =\dfrac{1}{2}\sum_{m,n=-\infty}^{\infty}(m^{2}-4n^{2})q^{n^{2}+(m^{2}-1)/4},\label{eq:c2-1}\\
(q;q)_{\infty}^{6} & =\dfrac{1}{2}\sum_{m,n=-\infty}^{\infty}(-1)^{m}(4n^{2}-m^{2})q^{n^{2}+(m^{2}-1)/4}.\nonumber 
\end{align}
\end{prop}
\noindent{\it Proof.} It can seen from the series expansions of $\theta_{2}$
and $\theta_{3}$ that 
\begin{equation}
\begin{gathered}\vartheta_{2}(2\tau)=q^{\frac{1}{2}}\sum_{k=-\infty}^{\infty}q^{2k(k+1)},\;\vartheta_{3}(\tau/2)=\sum_{k=-\infty}^{\infty}q^{k^{2}/2},\\
\vartheta''_{2}(2\tau)=-q^{\frac{1}{2}}\sum_{k=-\infty}^{\infty}(2k+1)^{2}q^{2k(k+1)},\;\vartheta''_{3}(\tau/2)=-4\sum_{k=-\infty}^{\infty}k^{2}q^{k^{2}/2}.
\end{gathered}
\label{eq:7-9}
\end{equation}
Differentiating both sides of (\ref{eq:1-25}) with respect to $x$
twice and then putting $x=y=0$ we have 
\[
2(\vartheta'_{1}(\tau))^{2}=\vartheta_{2}(2\tau)\vartheta''_{3}(\tau/2)-4\vartheta_{3}(\tau/2)\vartheta''_{2}(2\tau).
\]
Substituting (\ref{eq:7-3}) and (\ref{eq:7-9}) into this equation
and replacing $q^{2}$ by $q$ we easily obtain the first result.

Similarly, from (\ref{eq:1-26})\textendash (\ref{eq:1-28}) we can
also deduce the other three results. This concludes the proof of Corollary
\ref{c2}. \qed

We now use (\ref{eq:c2-1}) to prove the well-known Ramanujan partition
congruence\footnote{Actually, this congruence can be proved by any of the series expansions
for $(q;q)_{\infty}^{6}$ in Corollary \ref{c2}.}: 
\begin{equation}
p(7n+5)\equiv0(\bmod7),\label{eq:7-10}
\end{equation}
where $p(n)$ denotes the number of unrestricted partitions of the
positive integers $n.$

\noindent{\it Proof of \eqref{eq:7-10}.} Set $(q;q)_{\infty}^{6}=\sum_{n=0}^{\infty}a(n)q^{n}.$
Comparing the coefficients of $q^{n}$ on both sides of (\ref{eq:c2-1}),
we get 
\[
a(n)=\dfrac{1}{2}\sum_{\substack{u,v=-\infty\\
4v^{2}+u^{2}-1=4n
}
}^{\infty}(u^{2}-4v^{2}).
\]
If $n\equiv5(\bmod7),$ then $4v^{2}+u^{2}\equiv0(\bmod7).$ Examining
all cases modulo $7,$ we find that $u\equiv v\equiv0(\bmod7)$ and
so $u^{2}-4v^{2}\equiv0(\bmod7^{2}).$ This means that 
\begin{equation}
a(7n+5)\equiv0(\bmod7^{2}).\label{eq:7-11}
\end{equation}
By the binomial theorem, 
\begin{equation}
(q;q)_{\infty}^{7}\equiv(q^{7};q^{7})_{\infty}(\bmod7).\label{eq:7-12}
\end{equation}
From \cite[(1.1.7)]{B} and (\ref{eq:7-12}) we deduce that 
\[
\sum_{n=0}^{\infty}p(n)q^{n}=\dfrac{1}{(q;q)_{\infty}}=\dfrac{(q;q)_{\infty}^{6}}{(q;q)_{\infty}^{7}}\equiv\dfrac{\sum_{n=0}^{\infty}a(n)q^{n}}{(q^{7};q^{7})_{\infty}}(\bmod7).
\]
Extracting the terms with indices of the form $7n+5$ from the above
identity and using (\ref{eq:7-11}), we obtain 
\[
\sum_{n=0}^{\infty}p(7n+5)q^{7n+5}\equiv\dfrac{\sum_{n=0}^{\infty}a(7n+5)q^{7n+5}}{(q^{7};q^{7})_{\infty}}\equiv0(\bmod7).
\]
So $p(7n+5)\equiv0(\bmod7).$ This proves (\ref{eq:7-10}). \qed

\section{\label{sec:New} Proof of Theorem \ref{t2-1} and a series representation
for $\eta^{6}(q)$ }

We first prove Theorem $\ref{t2-1}$ using Theorem $\ref{t1-2}.$ 

\noindent{\it Proof of Theorem \ref{t2-1}.} It follows from $\eqref{eq:1-2}$
and $\eqref{eq:1-3}$ that $\theta_{3}(2z|2\tau)$ and $\theta_{2}(2z|2\tau)$
satisfy the conditions of Theorem $\ref{t1-2}.$ Putting $h_{1}(z)=\theta_{3}(2z|2\tau)$
and $h_{2}(z)=\theta_{2}(2z|2\tau)$ in Theorem $\ref{t1-2}$ gives
\begin{equation}
4\theta_{3}(2x|2\tau)\theta_{2}(2y|2\tau)-4\theta_{2}(2x|2\tau)\theta_{3}(2y|2\tau)=C\theta_{1}(x-y|\tau)\theta_{1}(x+y|\tau).\label{eq:2-3}
\end{equation}
We now determine the constant $C.$ Setting $x=\frac{\pi}{2},\;y=0$
in $\eqref{eq:2-3}$ and applying $\eqref{eq:1-2}$,$\eqref{eq:1-3}$
and $\eqref{eq:1-8}$, we arrive at 
\begin{equation}
8\vartheta_{2}(2\tau)\vartheta_{3}(2\tau)=C\vartheta_{2}^{2}(\tau).\label{eq:22-5}
\end{equation}
Applying the infinite product representations of theta functions,
we can deduce that
\[
\vartheta_{2}(\tau)\theta_{2}(2z|\tau)=2\theta_{2}(2z|2\tau)\theta_{3}(2z|2\tau).
\]
Setting $z=0$ in the above identity gives
\begin{equation}
\vartheta_{2}^{2}(\tau)=2\vartheta_{2}(2\tau)\vartheta_{3}(2\tau).\label{eq:1-11}
\end{equation}
We use $\eqref{eq:1-11}$ in (\ref{eq:22-5}) to get 
\[
C=4.
\]
Substituting $C=4$ into $\eqref{eq:2-3}$ and making the substitutions
$x+y\rightarrow x,\;x-y\rightarrow y,$ we readily obtain $\eqref{eq:2-1}$.
This finishes the proof of Theorem $\ref{t2-1}.$ \qed

We next derive a series representation for $\eta^{6}(q)$ using (\ref{eq:a-4}).
\begin{prop}
\label{c3} For $|q|<1,$ we have
\[
(q;q)_{\infty}^{6}=\dfrac{1}{2}\sum_{m,n=-\infty}^{\infty}((2n+1)^{2}-(2m)^{2})q^{m^{2}+n(n+1)}.
\]
\end{prop}
\noindent{\it Proof.} Differentiating both sides of (\ref{eq:a-4})
with respect to $x$ twice and then putting $x=y=0$ we have
\begin{equation}
(\vartheta'_{1}(\tau))^{2}=2\vartheta_{2}(2\tau)\vartheta''_{3}(2\tau)-2\vartheta{}_{3}(2\tau)\vartheta''_{2}(2\tau).\label{eq:a-2}
\end{equation}
It follows from the definition of the theta function $\theta_{3}$
that 
\[
\vartheta_{3}(2\tau)=\sum_{k=-\infty}^{\infty}q^{2k^{2}},\quad\vartheta''_{3}(2\tau)=-\sum_{k=-\infty}^{\infty}(2k)^{2}q^{2k^{2}}.
\]
Substituting these two identities, \eqref{eq:7-7} and the first and
the third formulas of \eqref{eq:7-9} into (\ref{eq:a-2}) and replacing
$q^{2}$ by $q$ we get
\[
(q;q)_{\infty}^{6}=\dfrac{1}{2}\sum_{m,n=-\infty}^{\infty}((2n+1)^{2}-(2m)^{2})q^{m^{2}+n(n+1)}.
\]
 This ends the proof of Proposition \ref{c3}. \qed 

\section{\label{sec:sn} New proofs for Jacobi's two, four and eight squares
theorems}

Representing natural numbers as sums of squares is an important topic
in number theory. Given a natural number $n$, denote by $r_{l}(n)$
the number of integer solutions of the Diophantine equation 
\[
n=x_{1}^{2}+x_{2}^{2}+x_{3}^{2}+\cdots+x_{l}^{2},
\]
which counts the number of ways in which $n$ can be written as a
sum of $l$ squares. In $l$-dimensional space, $r_{l}(n)$ also gives
the number of points with integer coordinates on the sphere. Jacobi
\cite{J} was the first to give the formulas for $r_{2}(n),\;r_{4}(n),\;r_{6}(n)$
and $r_{8}(n)$. From then on, many proofs for these formulas appeared,
see, for example, Andrews, Ekhad and Zeilberger \cite{AEZ}, Bhargava
and Adiga \cite{BA}, Chan \cite{Ch}, Cooper and Lam \cite{CL},
Hardy and Wright \cite[pp. 241--242]{HW}, Hirschhorn \cite{Hi1,Hi2},
Lin \cite{Lin} and Spearman and Williams \cite{SW}. In this section
we give new proofs for Jacobi's two, four and eight squares theorems
using $\eqref{eq:2-1}$ and $\eqref{eq:2-2}$. 
\begin{thm}
(Jacobi's Two Squares Theorem) For each positive integer $n$, we
have 
\[
r_{2}(n)=4\bigg(\sum_{\substack{d|n\\
d\equiv1(\bmod4)
}
}1-\sum_{\substack{d|n\\
d\equiv3(\bmod4)
}
}1\bigg).
\]
\end{thm}

\noindent{\it Proof.} Differentiating both sides of $\eqref{eq:2-1}$
with respect to $x$ and then setting $x=0,$ we find 
\begin{equation}
\begin{aligned}\vartheta'_{1}(\tau)\theta_{1}(y|\tau) & =\bigg(\frac{\theta'_{3}}{\theta_{3}}(y|2\tau)+\frac{\theta'_{2}}{\theta_{2}}(-y|2\tau)\bigg)\theta_{2}(-y|2\tau)\theta_{3}(y|2\tau)\\
 & \quad-\bigg(\frac{\theta'_{3}}{\theta_{3}}(-y|2\tau)+\frac{\theta'_{2}}{\theta_{2}}(y|2\tau)\bigg)\theta_{2}(y|2\tau)\theta_{3}(-y|2\tau)\\
 & =2\bigg(\frac{\theta'_{3}}{\theta_{3}}(y|2\tau)-\frac{\theta'_{2}}{\theta_{2}}(y|2\tau)\bigg)\theta_{2}(y|2\tau)\theta_{3}(y|2\tau).
\end{aligned}
\label{eq:3-3}
\end{equation}
Using the infinite product representations of theta functions, $\eqref{eq:1-5-1}$
and $\eqref{eq:1-6}$ in $\eqref{eq:3-3}$ and after some direct computations,
we get 
\begin{equation}
\begin{aligned} & \frac{\tan y(q^{2};q^{2})_{\infty}^{4}(q^{2}e^{2iy};q^{2})_{\infty}(q^{2}e^{-2iy};q^{2})_{\infty}}{(q^{4};q^{4})_{\infty}^{2}(-q^{2}e^{2iy};q^{2})_{\infty}(-q^{2}e^{-2iy};q^{2})_{\infty}}\\
 & =\tan y+4\sum_{n=1}^{\infty}\frac{(-1)^{n}q^{2n}}{1+q^{2n}}\sin(2ny).
\end{aligned}
\label{eq:3-7}
\end{equation}
It follows from the infinite product representations for $\vartheta_{4}(\tau)$
that 
\begin{equation}
\vartheta_{4}(\tau)=\frac{(q;q)_{\infty}}{(-q;q)_{\infty}}=\sum_{n=-\infty}^{\infty}(-1)^{n}q^{n^{2}}.\label{eq:1-7}
\end{equation}
Setting $y=\frac{\pi}{4}$ in $\eqref{eq:3-7}$, replacing $q^{2}$
by $q$ and applying $\eqref{eq:1-7}$ we are led to 
\[
\bigg(\sum_{n=-\infty}^{\infty}(-1)^{n}q^{n^{2}}\bigg)^{2}=1+4\sum_{n=1}^{\infty}\bigg(\frac{q^{4n-1}}{1+q^{4n-1}}-\frac{q^{4n-3}}{1+q^{4n-3}}\bigg).
\]
Replacing $q$ by $-q$ in this equation we arrive at 
\[
\begin{alignedat}{1}\bigg(\sum_{n=-\infty}^{\infty}q^{n^{2}}\bigg)^{2} & =1-4\sum_{n=1}^{\infty}\bigg(\frac{q^{4n-1}}{1-q^{4n-1}}-\frac{q^{4n-3}}{1-q^{4n-3}}\bigg)\\
 & =1-4\sum_{n=1}^{\infty}\sum_{k=1}^{\infty}\Big(q^{k(4n-1)}-q^{k(4n-3)}\Big)\\
 & =1+4\sum_{n=1}^{\infty}\bigg(\sum_{\substack{d|n\\
d\equiv1(\bmod4)
}
}1-\sum_{\substack{d|n\\
d\equiv3(\bmod4)
}
}1\bigg)q^{n}.
\end{alignedat}
\]
Comparing the coefficient of $q^{n}$ in this identity we readily
obtain the result. \qed
\begin{thm}
(Jacobi's Four Squares Theorem) For each positive integer $n$, we
have 
\[
r_{4}(n)=8\sum_{\substack{d|n\\
4\nmid d
}
}d.
\]
\end{thm}

\noindent{\it Proof.} It follows from $\eqref{eq:1-5-1}$ and $\eqref{eq:1-6}$
that
\begin{align}
\frac{\theta'_{3}}{\theta_{3}}(0|\tau) & =\frac{\theta'_{2}}{\theta_{2}}(0|\tau)=0,\label{eq:3-4}\\
\bigg(\frac{\theta'_{2}}{\theta_{2}}\bigg)'(0|\tau) & =-1+8\sum_{n=1}^{\infty}(-1)^{n}\frac{nq^{2n}}{1-q^{2n}},\label{eq:3-5}\\
\bigg(\frac{\theta'_{3}}{\theta_{3}}\bigg)'(0|\tau) & =8\sum_{n=1}^{\infty}(-1)^{n}\frac{nq^{n}}{1-q^{2n}}.\label{eq:3-6}
\end{align}
Differentiating both sides of $\eqref{eq:3-3}$ with respect to $y$
using the method of logarithmic differentiation, setting $y=0$ and
using $\eqref{eq:3-4}$ we get 
\begin{equation}
\begin{aligned}\text{(}\vartheta'_{1}(\tau))^{2} & =2\bigg(\bigg(\frac{\theta'_{3}}{\theta_{3}}\bigg)'(0|2\tau)-\bigg(\frac{\theta'_{2}}{\theta_{2}}\bigg)'(0|2\tau)\bigg)\vartheta_{2}(2\tau)\vartheta_{3}(2\tau)\\
 & \quad+2\bigg(\bigg(\frac{\theta'_{3}}{\theta_{3}}(0|2\tau)\bigg)^{2}-\bigg(\frac{\theta'_{2}}{\theta_{2}}(0|2\tau)\bigg)^{2}\bigg)\vartheta_{2}(2\tau)\vartheta_{3}(2\tau)\\
 & =2\bigg(\bigg(\frac{\theta'_{3}}{\theta_{3}}\bigg)'(0|2\tau)-\bigg(\frac{\theta'_{2}}{\theta_{2}}\bigg)'(0|2\tau)\bigg)\vartheta_{2}(2\tau)\vartheta_{3}(2\tau).
\end{aligned}
\label{eq:33-1}
\end{equation}
Using the infinite product representations of theta functions, we
can verify that 
\begin{eqnarray*}
 &  & \vartheta'_{1}(\tau)\theta_{1}(2z|\tau)=2\theta_{1}(z|\tau)\theta_{2}(z|\tau)\theta_{3}(z|\tau)\theta_{4}(z|\tau),\\
 &  & \theta_{3}(z|\tau)\theta_{4}(z|\tau)=\vartheta_{4}(2\tau)\theta_{4}(2z|2\tau).
\end{eqnarray*}
 Then 
\begin{align}
\vartheta'_{1}(\tau) & =\vartheta_{2}(\tau)\vartheta_{3}(\tau)\vartheta_{4}(\tau),\label{eq:1-12}\\
\vartheta_{4}^{2}(2\tau) & =\vartheta_{3}(\tau)\vartheta_{4}(\tau),\label{eq:1-13}
\end{align}
where (\ref{eq:1-12}) follows by differentiating the first identity
with respect to $z$ and then setting $z=0.$ Applying $\eqref{eq:1-11}$,$\eqref{eq:1-12}$
and $\eqref{eq:1-13}$ in (\ref{eq:33-1}), replacing $q^{2}$ by
$q$ and using $\eqref{eq:3-5}$ and $\eqref{eq:3-6}$ we have 
\[
\vartheta_{4}^{4}(\tau)=1+8\sum_{n=1}^{\infty}(-1)^{n}\frac{nq^{n}}{1+q^{n}}.
\]
Replacing $q$ by $-q$ in the above identity and using $\eqref{eq:1-7}$,
we obtain 
\begin{align*}
\bigg(\sum_{n=-\infty}^{\infty}q^{n^{2}}\bigg)^{4} & =1+8\sum_{n=1}^{\infty}\frac{nq^{n}}{1+(-q)^{n}}\\
 & =1+8\sum_{n=1}^{\infty}\frac{(2n-1)q^{2n-1}}{1-q^{2n-1}}+8\sum_{n=1}^{\infty}\frac{2nq^{2n}}{1+q^{2n}}\\
 & =1+8\sum_{n=1}^{\infty}\frac{(2n-1)q^{2n-1}}{1-q^{2n-1}}+8\sum_{n=1}^{\infty}\frac{2nq^{2n}}{1-q^{2n}}\\
 & \;+8\sum_{n=1}^{\infty}\frac{2nq^{2n}}{1+q^{2n}}-8\sum_{n=1}^{\infty}\frac{2nq^{2n}}{1-q^{2n}}\\
 & =1+8\sum_{n=1}^{\infty}\frac{nq^{n}}{1-q^{n}}-32\sum_{n=1}^{\infty}\frac{nq^{4n}}{1-q^{4n}}\\
 & =1+8\sum_{n=1}^{\infty}(\sum_{d|n}d-\sum_{\substack{d|n\\
4|d
}
}d)q^{n}.
\end{align*}
Equating the coefficient of $q^{n}$ we establish the result. \qed
\begin{thm}
(Jacobi's Eight Squares Theorem) For each positive integer $n$, we
have 
\[
r_{8}(n)=16\sum_{\substack{d|n}
}(-1)^{n+d}d^{3}.
\]
\end{thm}

\noindent{\it Proof.} Applying the infinite product representations
of theta functions, we can deduce that 
\begin{equation}
\theta_{1}(z|\tau)\theta_{2}(z|\tau)=\vartheta_{4}(2\tau)\theta_{1}(2z|2\tau),\label{eq:1-4}
\end{equation}
Multiplying $\eqref{eq:2-1}$ by $\eqref{eq:2-2}$, applying $\eqref{eq:1-4}$
and replacing $q^{2}$ by $q$, we get 
\[
\vartheta_{4}^{2}(\tau)\theta_{1}(2x|\tau)\theta_{1}(2y|\tau)=\theta_{2}^{2}(x-y|\tau)\theta_{3}^{2}(x+y|\tau)-\theta_{2}^{2}(x+y|\tau)\theta_{3}^{2}(x-y|\tau).
\]
Making the substitutions $x\rightarrow x+\frac{\pi\tau}{4},\;y\rightarrow y+\frac{\pi\tau}{4}$
in the above indetity and using $\eqref{eq:1-8}$\textendash $\eqref{eq:1-10}$
we have 
\[
-\vartheta_{4}^{2}(\tau)\theta_{4}(2x|\tau)\theta_{4}(2y|\tau)=\theta_{2}^{2}(x-y|\tau)\theta_{2}^{2}(x+y|\tau)-\theta_{3}^{2}(x-y|\tau)\theta_{3}^{2}(x+y|\tau).
\]
Replacing $x$ by $x+\frac{\pi}{2}$ in the above identity, applying
$\eqref{eq:1-14}$, $\eqref{eq:1-9}$, $\eqref{eq:1-10}$ and $\eqref{eq:1-7}$
and setting $y=0$ we deduce 
\[
-\theta_{4}(2x|\tau)\vartheta_{4}^{3}(\tau)=\theta_{1}^{4}(x|\tau)-\theta_{4}^{4}(x|\tau).
\]
Then 
\begin{equation}
-\theta_{4}\bigg(\frac{2x}{\tau}\bigg|-\frac{1}{\tau}\bigg)\vartheta_{4}^{3}\bigg(-\frac{1}{\tau}\bigg)=\theta_{1}^{4}\bigg(\frac{x}{\tau}\bigg|-\frac{1}{\tau}\bigg)-\theta_{4}^{4}\bigg(\frac{x}{\tau}\bigg|-\frac{1}{\tau}\bigg).\label{eq:3-1}
\end{equation}
Recall the imaginary transformation formulas (see \cite[p. 475]{WW}
and \cite[p. 140]{L4}): 
\begin{align*}
 & \theta_{1}\bigg(\frac{z}{\tau}\bigg|-\frac{1}{\tau}\bigg)=-i\sqrt{-i\tau}\exp((iz^{2})/(\pi\tau))\theta_{1}(z|\tau),\\
 & \theta_{4}\bigg(\frac{z}{\tau}\bigg|-\frac{1}{\tau}\bigg)=\sqrt{-i\tau}\exp((iz^{2})/(\pi\tau))\theta_{2}(z|\tau).
\end{align*}
Applying these imaginary transformation formulas in $\eqref{eq:3-1}$
and cancelling out certain common factors gives 
\[
\theta_{2}(2x|\tau)\vartheta_{2}^{3}(\tau)=\theta_{2}^{4}(x|\tau)-\theta_{1}^{4}(x|\tau).
\]
Differentiating this identity four times with respect to $x$ using
the method of logarithmic differentiation and then setting $x=0$,
we find 
\[
\bigg(\frac{\theta'_{2}}{\theta_{2}}\bigg)'''(0|\tau)=-\frac{2[\vartheta'_{1}(\tau)]^{4}}{\vartheta_{2}^{4}(\tau)}.
\]
Using $\eqref{eq:1-12}$, $\eqref{eq:1-13}$ and $\eqref{eq:1-5-1}$
in the above identity and replacing $q^{2}$ by $q$, we arrive at
\[
\vartheta_{4}^{8}(\tau)=1+16\sum_{n=1}^{\infty}\frac{n^{3}(-q)^{n}}{1-q^{n}}.
\]
We use $\eqref{eq:1-7}$ and replace $-q$ by $q$ in the above identity
to get 
\begin{equation}
\begin{aligned}\bigg(\sum_{n=-\infty}^{\infty}q^{n^{2}}\bigg)^{8} & =1+16\sum_{n=1}^{\infty}\frac{n^{3}q{}^{n}}{1-(-q)^{n}}\\
 & =1+16\sum_{n=1}^{\infty}\frac{(2n-1)^{3}q{}^{2n-1}}{1+q^{2n-1}}+16\sum_{n=1}^{\infty}\frac{(2n)^{3}q{}^{2n}}{1-q^{2n}}\\
 & =1+16\sum_{n=1}^{\infty}\frac{(2n)^{3}q{}^{2n}}{1-q^{2n}}+16\sum_{n=1}^{\infty}\frac{(2n-1)^{3}q{}^{2n-1}}{1-q^{2n-1}}\\
 & \quad+16\sum_{n=1}^{\infty}\frac{(2n-1)^{3}q{}^{2n-1}}{1+q^{2n-1}}-16\sum_{n=1}^{\infty}\frac{(2n-1)^{3}q{}^{2n-1}}{1-q^{2n-1}}\\
 & =1+16\sum_{n=1}^{\infty}\frac{n^{3}q{}^{n}}{1-q^{n}}-32\sum_{n=1}^{\infty}\frac{(2n-1)^{3}q{}^{4n-2}}{1-q^{4n-2}}\\
 & =1+16\sum_{n=1}^{\infty}\Big(\sum_{\substack{d|n}
}d^{3}-2\sum_{\substack{2d|n\\
d\;odd
}
}d^{3}\Big)q^{n}.
\end{aligned}
\label{eq:3-2}
\end{equation}
When $n$ is even,

\begin{align*}
\sum_{\substack{d|n}
}d^{3}-2\sum_{\substack{2d|n\\
d~odd
}
}d^{3} & =\sum_{\substack{d|n}
}d^{3}-2\sum_{\substack{d|n\\
d~odd
}
}d^{3}\\
 & =\sum_{\substack{d|n\\
d~even
}
}d^{3}+\sum_{\substack{d|n\\
d~odd
}
}d^{3}-2\sum_{\substack{d|n\\
d~odd
}
}d^{3}\\
 & =\sum_{\substack{d|n}
}(-1)^{d}d^{3}.
\end{align*}
When $n$ is odd, 
\[
\sum_{\substack{d|n}
}d^{3}-2\sum_{\substack{2d|n\\
d\;odd
}
}d^{3}=\sum_{\substack{d|n}
}d^{3}.
\]
Therefore, for any positive integer $n$ we have

\[
\sum_{\substack{d|n}
}d^{3}-2\sum_{\substack{2d|n\\
d\;odd
}
}d^{3}=\sum_{\substack{d|n}
}(-1)^{n+d}d^{3}.
\]
Then the result follows readily by comparing the coefficient of $q^{n}$
in $\eqref{eq:3-2}$. \qed

We now use $\eqref{eq:3-7}$ to give a new proof of a well-known result.
Our proof is different from that of Berndt \cite{B}.
\begin{thm}
Let $r_{1,3}(n)$ denote the number of representations of the positive
integer $n$ by $x^{2}+3y^{2}.$ Then
\[
r_{1,3}(n)=4\bigg(\sum_{\substack{d|n\\
d\equiv4(\bmod12)
}
}1-\sum_{\substack{d|n\\
d\equiv8(\bmod12)
}
}1\bigg)+2\bigg(\sum_{\substack{d|n\\
d\equiv1(\bmod3)
}
}1-\sum_{\substack{d|n\\
d\equiv2(\bmod3)
}
}1\bigg).
\]
\end{thm}

\noindent{\it Proof.}  It follows from the Jacobi triple product
identity that 
\[
\sum_{n=-\infty}^{\infty}(-1)^{n}q^{n^{2}}=(q^{2};q^{2})_{\infty}(q;q^{2})_{\infty}^{2}=\dfrac{(q^{2};q^{2})_{\infty}}{(-q;q)_{\infty}^{2}}=\dfrac{(q;q)_{\infty}}{(-q;q)_{\infty}}.
\]
Then
\begin{equation}
\sum_{n=-\infty}^{\infty}q^{n^{2}}=\dfrac{(-q;-q)_{\infty}}{(q;-q)_{\infty}}.\label{eq:3-10}
\end{equation}
Setting $y=\frac{\pi}{3}$ in $\eqref{eq:3-7}$, replacing $q^{2}$
by $q$ and making the substitution $q\rightarrow-q$ we derive 
\[
\dfrac{(-q;q)_{\infty}(-q^{3};-q^{3})_{\infty}}{(q;-q)_{\infty}(q^{3};-q^{3})_{\infty}}=1+2\sum_{n=1}^{\infty}\bigg(\frac{n}{3}\bigg)\frac{q^{n}}{1+(-q)^{n}},
\]
where $\big(\frac{.}{3}\big)$ denotes the Legendre symbol. Applying
(\ref{eq:3-10}) in the above identity, we have
\begin{align*}
\sum_{n=-\infty}^{\infty}q^{n^{2}}\cdot\sum_{n=-\infty}^{\infty}q^{3n^{2}} & =1+2\sum_{n=1}^{\infty}\bigg(\frac{2n}{3}\bigg)\frac{q^{2n}}{1+q^{2n}}+2\sum_{n=0}^{\infty}\bigg(\frac{2n+1}{3}\bigg)\frac{q^{2n+1}}{1-q^{2n+1}}\\
 & =1+2\sum_{n=1}^{\infty}\bigg(\frac{2n}{3}\bigg)\frac{q^{2n}}{1-q^{2n}}+2\sum_{n=0}^{\infty}\bigg(\frac{2n+1}{3}\bigg)\frac{q^{2n+1}}{1-q^{2n+1}}\\
 & \quad-4\sum_{n=1}^{\infty}\bigg(\frac{2n}{3}\bigg)\frac{q^{4n}}{1-q^{4n}}\\
 & =1+2\sum_{n=1}^{\infty}\bigg(\frac{n}{3}\bigg)\frac{q^{n}}{1-q^{n}}+4\sum_{n=1}^{\infty}\bigg(\frac{n}{3}\bigg)\frac{q^{4n}}{1-q^{4n}}\\
 & =1+2\sum_{n=0}^{\infty}\bigg(\frac{q^{3n+1}}{1-q^{3n+1}}-\frac{q^{3n+2}}{1-q^{3n+2}}\bigg)\\
 & \quad+4\sum_{n=0}^{\infty}\bigg(\frac{q^{4(3n+1)}}{1-q^{4(3n+1)}}-\frac{q^{4(3n+2)}}{1-q^{4(3n+2)}}\bigg).
\end{align*}
Expanding each of the quotients $1/(1-q^{k})$ into geometric series
and then equating the coefficients of $q^{n}$ on both sides of the
above identity we can easily obtain the result. \qed

\section{\label{sec:Proof1} Proofs of Theorems $\ref{t2-2}$\textendash \ref{t2-4}
and $q$-Trigonometric identities}

In this section we first prove Theorems $\ref{t2-2}$, $\ref{t2-3}$
and \ref{t2-4} using Theorem $\ref{t1-2}$.

\noindent{\it Proof of Theorem \ref{t2-2}.} From $\eqref{eq:1-1}$
and $\eqref{eq:1-2}$ we see that $\theta_{1}(z-u|\tau)\theta_{1}(z+u|\tau)$
and $\theta_{2}(z-v|\tau)\theta_{2}(z+v|\tau)$ satisfy the conditions
of Theorem $\ref{t1-2}$. Taking $h_{1}(z)=\theta_{1}(z-u|\tau)\theta_{1}(z+u|\tau)$
and $h_{2}(z)=\theta_{2}(z-v|\tau)\theta_{2}(z+v|\tau)$ in Theorem
$\ref{t1-2}$, we get 
\begin{align*}
 & 4\theta_{1}(x-u|\tau)\theta_{1}(x+u|\tau)\theta_{2}(y-v|\tau)\theta_{2}(y+v|\tau)\\
 & \quad-4\theta_{1}(y-u|\tau)\theta_{1}(y+u|\tau)\theta_{2}(x-v|\tau)\theta_{2}(x+v|\tau)\\
 & =C\theta_{1}(x-y|\tau)\theta_{1}(x+y|\tau).
\end{align*}
We put $y=u$ in the above identity to know that $C=4\theta_{2}(u-v|\tau)\theta_{2}(u+v|\tau)$
and then substitute this constant back into the equality to obtain
(\ref{eq:2.1}). \qed

\noindent{\it Proof of Theorem \ref{t2-3}.} We first show $\eqref{eq:22-4}$.
Using $\eqref{eq:1-1}$ and $\eqref{eq:1-2}$ we can verify $\theta_{1}(z-u|\tau)\theta_{1}(z+u|\tau)$
and $\theta_{2}(2z|2\tau)$ satisfy the functional equations $\eqref{eq:t-1}$.
Setting $h_{1}(z)=\theta_{1}(z-u|\tau)\theta_{1}(z+u|\tau)$ and $h_{2}(z)=\theta_{2}(2z|2\tau)$
in Theorem $\ref{t1-2}$, we have 
\begin{gather}
\begin{gathered}4\theta_{1}(x-u|\tau)\theta_{1}(x+u|\tau)\theta_{2}(2y|2\tau)-4\theta_{1}(y-u|\tau)\theta_{1}(y+u|\tau)\theta_{2}(2x|2\tau)\\
=C\theta_{1}(x-y|\tau)\theta_{1}(x+y|\tau).
\end{gathered}
\label{eq:4-1}
\end{gather}
Taking $y=u$ in the above identity we get 
\[
C=4\theta_{2}(2u|2\tau).
\]
We substitute $C=4\theta_{2}(2u|2\tau)$ into $\eqref{eq:4-1}$ to
obtain the result.

Similarly, we can also deduce $\eqref{eq:22-8}$. This completes the
proof of Theorem $\ref{t2-3}$. \qed

\noindent{\it Proof of Theorem \ref{t2-4}.} It follows from (\ref{eq:1-1})
that the entire functions $\frac{\theta_{1}(z|\tau/3)}{\theta_{1}(z|\tau)}$
and $\theta_{1}^{2}(z|\tau)$ satisfy the functional equations (\ref{eq:t-1}).
Setting $h_{1}(z)=\frac{\theta_{1}(z|\tau/3)}{\theta_{1}(z|\tau)}$
and $h_{2}(z)=\theta_{1}^{2}(z|\tau)$ in Theorem \ref{t1-2} we get
\begin{equation}
\frac{4\theta_{1}(x|\tau/3)\theta_{1}^{2}(y|\tau)}{\theta_{1}(x|\tau)}-\frac{4\theta_{1}^{2}(x|\tau)\theta_{1}(y|\tau/3)}{\theta_{1}(y|\tau)}=C\theta_{1}(x-y|\tau)\theta_{1}(x+y|\tau).\label{eq:7-1}
\end{equation}
Putting $y=0$ in this identity and then cancelling out the factor
$\theta_{1}^{2}(x|\tau)$ gives
\[
C=-4\lim_{y\rightarrow0}\frac{\theta_{1}(y|\tau/3)}{\theta_{1}(y|\tau)}=-4\frac{\vartheta'_{1}(\tau/3)}{\vartheta'_{1}(\tau)}.
\]
Substituting this constant back into (\ref{eq:7-1}) and simplifying
we can easily obtain (\ref{eq:1-20}). \qed

Gosper in \cite{G} introduced $q$-analogues of $\sin z$ and $\cos z$
which are defined as 
\begin{eqnarray}
 &  & \sin_{q}(\pi z):=q^{(z-1\text{/2})^{2}}\prod_{n=1}^{\infty}\frac{(1-q^{2n-2z})(1-q^{2n+2z-2})}{(1-q^{2n-1})^{2}},\nonumber \\
 &  & \cos_{q}(\pi z):=q^{z^{2}}\prod_{n=1}^{\infty}\frac{(1-q^{2n-2z-1})(1-q^{2n+2z-1})}{(1-q^{2n-1})^{2}}.\label{eq:4-7}
\end{eqnarray}
See \cite{GR} for other $q$-analogues of the trigonometric functions.
From the above definitions of $\sin_{q}$ and $\cos_{q},$ we can
see that $\cos_{q}z=\sin_{q}(\frac{\pi}{2}\pm z),\;\lim_{q\rightarrow1}\sin_{q}z=\sin z$
and $\lim_{q\rightarrow1}\cos_{q}z=\cos z.$ In \cite{G}, Gosper
gave two relations between $\sin_{q},\;\cos_{q}$ and the functions
$\theta_{1}$ and $\theta_{2},$ which are equivalent to the following
identities: 
\begin{equation}
\sin_{q}z=\frac{\theta_{1}(z|\tau')}{\vartheta_{2}(\tau')},\;\cos_{q}z=\frac{\theta_{2}(z|\tau')}{\vartheta_{2}(\tau')},\label{eq:4-2}
\end{equation}
where $\tau'=-\frac{1}{\tau}.$ In the same paper, He stated without
proofs various identities involving $\sin_{q}z$ and $\cos_{q}z.$
He conjectured these identities by using empirical evidence based
on a computer program called MACSYMA. Using a direct analysis of its
left-hand side and its right-hand side through logarithmic derivatives,
Mez\H{o} in \cite{M} confirmed the conjecture ($q$-Double$_{2}$).
Applying the theory of elliptic functions, M. El Bachraoui \cite{E}
proved the conjectures ($q$-Double$_{2}$) and ($q$-Double$_{3}$).
See \cite{AAE} for proofs of another two $q$-trigonometric identities
using theory of elliptic functions. Additionally, Gosper \cite[pp.100--103]{G}
conjectured the following identities on the $q$-trigonometric functions
and the constant $\Pi_{q}$:\footnote{Actually, in \cite[p. 101]{G} the $q$-trigonometric identity $\eqref{eq:4-6}$
is $\cos''_{q}0=(\sin'_{q}0)^{2}+\frac{\cos''_{\sqrt{q}}0}{2}.$ This
identity is not true since when $q\rightarrow1,$ it reduces to $\frac{\cos''0}{2}=(\sin'0)^{2},$
which obviouly does not hold. So we modify it here.} 
\begin{eqnarray}
 &  & \sin{}_{q}(x+y)\sin{}_{q}(x-y)=\sin{}_{q}^{2}x\cos_{q}^{2}y-\cos_{q}^{2}x\sin{}_{q}^{2}y,\label{eq:4-9}\\
 &  & (\sin'_{q}z)^{2}=\sin{}_{q}z\sin''_{q}z+(\sin'_{q}0)^{2}(\cos{}_{q}z)^{2}-\cos''_{q}0(\sin{}_{q}z)^{2},\label{eq:4-10}\\
 &  & \cos''_{q}0=(\sin'_{q}0)^{2}+2\cos''_{\sqrt{q}}0=\frac{2\ln q}{\pi^{2}}\bigg(1-4\ln q\sum_{n\geq1}\frac{q^{2n-1}}{(1-q^{2n-1})^{2}}\bigg),\label{eq:4-6}\\
 &  & \sqrt{\Pi_{q}\Pi_{q^{9}}}(\Pi_{q}+3\Pi_{q^{9}})=\Pi_{q^{3}}^{2}+3\Pi_{q}\Pi_{q^{9}},\label{eq:c2-4}
\end{eqnarray}
where $\Pi_{q}$ is given by 
\[
\Pi_{q}=q^{1/4}\frac{(q^{2};q^{2})_{\infty}^{2}}{(q;q^{2})_{\infty}^{2}}.
\]

In this section we shall establish three general $q$-trigonometric
identities which include as special cases $\eqref{eq:4-9}$, $\eqref{eq:4-10}$
and $\eqref{eq:4-6}$.

\begin{thm}
\label{t5-1} For all complex numbers $x,\;y,\;z,\;u$ and $v,$ we
have
\[
\begin{aligned} & \sin{}_{q}(x-u)\sin{}_{q}(x+u)\cos_{q}(y-v)\cos_{q}(y+v)\\
 & \quad-\cos_{q}(x-v)\cos_{q}(x+v)\sin{}_{q}(y-u)\sin{}_{q}(y+u)\\
 & =\cos_{q}(u-v)\cos_{q}(u+v)\sin{}_{q}(x+y)\sin{}_{q}(x-y)
\end{aligned}
\]
and 
\[
\begin{aligned} & \cos_{q}(z-v)\cos_{q}(z+v)\big((\sin'_{q}u)^{2}-\sin''_{q}u\sin{}_{q}u\big)\\
 & \quad+\sin{}_{q}(z-u)\sin{}_{q}(z+u)\big((\cos'_{q}v)^{2}-\cos''_{q}v\cos{}_{q}v\big)\\
 & =\cos_{q}(u-v)\cos_{q}(u+v)\big((\sin'_{q}z)^{2}-\sin''_{q}z\sin{}_{q}z\big).
\end{aligned}
\]
\end{thm}

\noindent{\it Proof.} Dividing both sides of the identity in Theorem
$\ref{t2-2}$ by $\vartheta_{2}^{4}(\tau),$ replacing $\tau$ by
$\tau'$ and applying $\eqref{eq:4-2}$ we obtain the first identity
immediately.

We now prove the second identity. Differentiating both sides of the
identity in Theorem $\ref{t2-2}$ twice with respect to $x$ and then
setting $x=0$ we get 
\begin{align*}
 & \theta_{2}(y-v|\tau)\theta_{2}(y+v|\tau)\big((\theta'_{1}(u|\tau))^{2}-\theta''_{1}(u|\tau)\theta_{1}(u|\tau)\big)\\
 & \quad+\theta_{1}(y-u|\tau)\theta_{1}(y+u|\tau)\big((\theta'_{2}(v|\tau))^{2}-\theta''_{2}(v|\tau)\theta_{2}(v|\tau)\big)\\
 & =\theta_{2}(u-v|\tau)\theta_{2}(u+v|\tau)\big((\theta'_{1}(y|\tau))^{2}-\theta''_{1}(y|\tau)\theta_{1}(y|\tau)\big).
\end{align*}
Dividing both sides of the above identity by $\vartheta_{2}^{4}(\tau)$,
making the substitution $\tau\rightarrow\tau'$, employing $\eqref{eq:4-2}$
and replacing $y$ by $z$ we achieve the result readily. This finishes
the proof of Theorem $\ref{t5-1}$. \qed
\begin{cor}
The $q$-trigonometric identities $\eqref{eq:4-9}$ and $\eqref{eq:4-10}$
are true.
\end{cor}
\noindent{\it Proof.} The $q$-trigonometric identities $\eqref{eq:4-9}$
and $\eqref{eq:4-10}$ follow easily from Theorem $\ref{t5-1}$ by
setting $u=v=0$ and applying the identities $\cos_{q}0=1$ and $\sin{}_{q}0=\cos'_{q}0=0$.
\qed

Replacing $x$ by $x+\frac{\pi}{2}$ in $\eqref{eq:4-9}$ we obtain
the result: 
\[
\cos{}_{q}(x+y)\cos{}_{q}(x-y)=\cos{}_{q}^{2}x\cos_{q}^{2}y-\sin_{q}^{2}x\sin{}_{q}^{2}y.
\]
When $y=x,$ this identity reduces to 
\[
\cos{}_{q}(2x)=\cos{}_{q}^{4}x-\sin_{q}^{4}x,
\]
which is also conjectured by Gosper \cite[p. 99]{G}.

When $q\rightarrow1,$ the first result in Theorem \ref{t5-1} reduces
to the following interesting trigonometric identity:
\[
\begin{aligned} & \sin(x-u)\sin(x+u)\cos(y-v)\cos(y+v)\\
 & \quad-\cos(x-v)\cos(x+v)\sin(y-u)\sin(y+u)\\
 & =\cos(u-v)\cos(u+v)\sin(x+y)\sin(x-y).
\end{aligned}
\]
\begin{thm}
\label{t4-3} For any complex number $z,$ we have 
\[
\cos''_{q}0\cos{}_{\sqrt{q}}(2z)=(\sin'_{q}z)^{2}-\sin{}_{q}z\sin_{q}''z+2\cos''_{\sqrt{q}}0\cos_{q}^{2}z.
\]
\end{thm}

\noindent{\it Proof.} It follows from $\eqref{eq:4-2}$ that 
\begin{equation}
\cos_{\sqrt{q}}z=\frac{\theta_{2}(z|2\tau')}{\vartheta_{2}(2\tau')}.\label{eq:4-5}
\end{equation}
Differentiating $\eqref{eq:22-4}$ twice with respect to $x$ and
then setting $x=0$ we have 
\begin{gather*}
2((\theta_{1}'(u|\tau))^{2}-\theta_{1}(u|\tau)\theta_{1}''(u|\tau))\theta_{2}(2y|2\tau)-4\vartheta_{2}''(2\tau)\theta_{1}(y-u|\tau)\theta_{1}(y+u|\tau)\\
=2\theta_{2}(2u|2\tau)((\theta'_{1}(y|\tau))^{2}-\theta_{1}(y|\tau)\theta''_{1}(y|\tau)).
\end{gather*}
We take $y=\frac{\pi}{2}$ in this identity and applying $\eqref{eq:1-2}$,
$\eqref{eq:1-8}$ and the identity $\theta'_{1}(\frac{\pi}{2}|\tau)=0$
to get 
\[
((\theta_{1}'(u|\tau))^{2}-\theta_{1}(u|\tau)\theta_{1}''(u|\tau))\vartheta_{2}(2\tau)+2\vartheta_{2}''(2\tau)\theta_{2}^{2}(u|\tau)=\theta_{2}(2u|2\tau)\vartheta_{2}(\tau)\vartheta''_{2}(\tau).
\]
Dividing both sides of the above identity by $\vartheta_{2}(2\tau)\vartheta_{2}^{2}(\tau)$,
making the substitution $\tau\rightarrow\tau'$, using $\eqref{eq:4-2}$
and $\eqref{eq:4-5}$ and replacing $u$ by $z$ we obtain 
\[
(\sin'_{q}z)^{2}-\sin{}_{q}z\sin_{q}''z+2\cos''_{\sqrt{q}}0\cos_{q}^{2}z=\cos''_{q}0\cos{}_{\sqrt{q}}(2z).
\]
This concludes the proof of Theorem $\ref{t4-3}$. \qed
\begin{cor}
The identity (\ref{eq:4-6}) is true.
\end{cor}
\noindent{\it Proof.} Setting $z=0$ in Theorem $\ref{t4-3}$ and
applying the identities $\sin{}_{q}0=0$ and $\cos_{q}0=1$ we can
deduce the first equality of $\eqref{eq:4-6}$ immediately. The second
equality of $\eqref{eq:4-6}$ follows readily by differentiating both
sides of $\eqref{eq:4-7}$ with respect to $z$ using the method of
logarithmic differentiation. \qed

We now use Theorem \ref{t2-4} to establish an interesting identity
involving the $q$-trigonometric functions and the constant $\Pi_{q}.$
From this identity we can deduce the identity (\ref{eq:c2-4}).
\begin{thm}
\label{t6-3} For all complex numbers $x$ and $y$ we have
\[
\mathrm{sin}_{q^{3}}y\mathrm{cos}_{q}^{3}x-\mathrm{cos}_{q^{3}}x\mathrm{sin}_{q}^{3}y=\frac{3\Pi_{q^{3}}}{\Pi_{q}}\mathrm{cos}_{q}x\mathrm{sin}_{q}y\mathrm{cos}_{q}(x-y)\mathrm{cos}_{q}(x+y).
\]
\end{thm}
In order to prove this theorem we need the following auxiliary result. 
\begin{lem}
\label{l1} For any positive integer $n,$ we have 
\[
\frac{\Pi_{q}}{\Pi_{q^{n}}}=\frac{n\vartheta'_{1}(\tau')\vartheta_{2}(\tau'/n)}{\vartheta'_{1}(\tau'/n)\vartheta_{2}(\tau')}.
\]
\end{lem}

\noindent{\it Proof.} Recall from \cite{WW} the following imaginary
transformation formula for the Jacobi theta function $\theta_{1}:$
\begin{align}
\theta_{1}(z\tau'|\tau') & =i\sqrt{-i\tau}\exp((iz^{2})/(\pi\tau))\theta_{1}(z|\tau),\label{eq:7}\\
\theta_{2}(z\tau'|\tau') & =\sqrt{-i\tau}\exp((iz^{2})/(\pi\tau))\theta_{4}(z|\tau),\label{eq:8}
\end{align}
Differentiating (\ref{eq:7}) with respect to $z$ and then setting
$z=0$ we have 
\[
\vartheta'_{1}(\tau')=-i\tau\sqrt{-i\tau}\vartheta'_{1}(\tau).
\]
Then
\[
\vartheta'_{1}(\tau'/n)=-in\tau\sqrt{-in\tau}\vartheta'_{1}(n\tau).
\]
So, by (\ref{eq:7-3}),
\[
\frac{\vartheta'_{1}(\tau')}{\vartheta'_{1}(\tau'/n)}=\frac{\vartheta'_{1}(\tau)}{n\sqrt{n}\vartheta'_{1}(n\tau)}=\frac{(q^{2};q^{2})_{\infty}^{3}}{n\sqrt{n}q^{\frac{n-1}{4}}(q^{2n};q^{2n})_{\infty}^{3}}.
\]

It follows from (\ref{eq:8}) that
\[
\vartheta_{2}(\tau')=\sqrt{-i\tau}\vartheta_{4}(\tau).
\]
Then 
\[
\vartheta_{2}(\tau'/n)=\sqrt{-in\tau}\vartheta_{4}(n\tau).
\]
Therefore, by (\ref{eq:7-4}), 
\[
\frac{\vartheta_{2}(\tau'/n)}{\vartheta_{2}(\tau')}=\frac{\sqrt{n}\vartheta_{4}(n\tau)}{\vartheta_{4}(\tau)}=\frac{\sqrt{n}(q^{2n};q^{2n})_{\infty}(q^{n};q^{2n})_{\infty}^{2}}{(q^{2};q^{2})_{\infty}(q;q^{2})_{\infty}^{2}}.
\]
In view of the above, we have
\[
\frac{\vartheta'_{1}(\tau')\vartheta_{2}(\tau'/n)}{\vartheta'_{1}(\tau'/n)\vartheta_{2}(\tau')}=\frac{(q^{2};q^{2})_{\infty}^{2}(q^{n};q^{2n})_{\infty}^{2}}{nq^{\frac{n-1}{4}}(q^{2n};q^{2n})_{\infty}^{2}(q;q^{2})_{\infty}^{2}}=\frac{\Pi_{q}}{n\Pi_{q^{n}}},
\]
which completes the proof of Lemma \ref{l1}. \qed

We are now in the position to prove Theorem \ref{t6-3}.

\noindent{\it Proof of Theorem \ref{t6-3}.} It is easily seen from
$\eqref{eq:4-2}$ that 
\begin{equation}
\mathrm{sin}_{q^{3}}z=\frac{\theta_{1}(y|\tau'/3)}{\vartheta_{2}(\tau'/3)},\:\mathrm{cos}_{q^{3}}z=\frac{\theta_{2}(y|\tau'/3)}{\vartheta_{2}(\tau'/3)}.\label{eq:77-1}
\end{equation}
Replacing $x$ by $x+\pi/2$ in (\ref{eq:1-20}) we have
\[
\begin{gathered}\theta_{2}^{3}(x|\tau)\theta_{1}(y|\tau/3)-\theta_{1}^{3}(y|\tau)\theta_{2}(x|\tau/3)\\
=\frac{\vartheta'_{1}(\tau/3)}{\vartheta'_{1}(\tau)}\theta_{2}(x|\tau)\theta_{1}(y|\tau)\theta_{2}(x-y|\tau)\theta_{2}(x+y|\tau).
\end{gathered}
\]
We replace $\tau$ by $\tau'$ in this equation, divide both sides
of the resulting identity by $\vartheta_{2}(\tau'/3)\vartheta_{2}^{3}(\tau')$
and then employ (\ref{eq:77-1}) to get 
\[
\mathrm{sin}_{q^{3}}y\mathrm{cos}_{q}^{3}x-\mathrm{cos}_{q^{3}}x\mathrm{sin}_{q}^{3}y=\frac{\vartheta'_{1}(\tau/3)\vartheta_{2}(\tau')}{\vartheta_{2}(\tau'/3)\vartheta'_{1}(\tau)}\mathrm{cos}_{q}x\mathrm{sin}_{q}y\mathrm{cos}_{q}(x-y)\mathrm{cos}_{q}(x+y).
\]
Then the result follows readily by combining the above identity and
the special case $n=3$ of Lemma \ref{l1}. \qed
\begin{cor}
The identity (\ref{eq:c2-4}) holds for $|q|<1$.
\end{cor}
\noindent{\it Proof.} We set $x=0$ in Theorem \ref{t6-3} to give
\begin{equation}
\mathrm{sin}_{q^{3}}y-\mathrm{sin}_{q}^{3}y=\frac{3\Pi_{q^{3}}}{\Pi_{q}}\mathrm{sin}_{q}y\mathrm{cos}_{q}^{2}y.\label{eq:77-2}
\end{equation}
According to \cite[p. 96]{G} we have
\begin{align}
\bigg(\mathrm{sin}_{q^{3}}\frac{\pi}{6}\bigg)^{3} & =\frac{1}{\big(\frac{\Pi_{q}}{\Pi_{q^{3}}}\big)^{2}-1},\label{eq:cc2-1}\\
\bigg(\mathrm{cos}_{q^{3}}\frac{\pi}{6}\bigg)^{3} & =\frac{\big(\frac{\Pi_{q}}{\Pi_{q^{3}}}\big)^{\frac{3}{2}}}{\big(\frac{\Pi_{q}}{\Pi_{q^{3}}}\big)^{2}-1}.\label{eq:c2-2}
\end{align}
We temporarily assume that $0<q<1$. Putting $z=\frac{\pi}{6}$ in
(\ref{eq:77-2}), replacing $q$ by $q^{3}$ and applying (\ref{eq:cc2-1})
and (\ref{eq:c2-2}) in the resulting identity we are led to
\begin{equation}
\frac{1}{\sqrt[3]{\big(\frac{\Pi_{q^{3}}}{\Pi_{q^{9}}}\big)^{2}-1}}-\frac{1}{\big(\frac{\Pi_{q}}{\Pi_{q^{3}}}\big)^{2}-1}=\frac{3\Pi_{q}\Pi_{q^{9}}}{\Pi_{q^{3}}^{2}}\frac{1}{\big(\frac{\Pi_{q}}{\Pi_{q^{3}}}\big)^{2}-1}.\label{eq:3-10-1}
\end{equation}
Let 
\[
x=\sqrt[3]{\bigg(\frac{\Pi_{q^{3}}}{\Pi_{q^{9}}}\bigg)^{2}-1},\;y=\sqrt{\frac{\Pi_{q}}{\Pi_{q^{9}}}}.
\]
Then
\[
x>-1,\;y>0,\;\frac{\Pi_{q^{3}}}{\Pi_{q^{9}}}=\sqrt{x^{3}+1}
\]
and
\[
\frac{\Pi_{q}}{\Pi_{q^{3}}}=\frac{\Pi_{q}}{\Pi_{q^{9}}}\bigg/\frac{\Pi_{q^{3}}}{\Pi_{q^{9}}}=\frac{y^{2}}{\sqrt{x^{3}+1}}.
\]
We substitute these identities into (\ref{eq:3-10-1}) and simplify
to get
\[
x=\frac{y^{4}-x^{3}-1}{x^{3}+3y^{2}+1}.
\]
This means that
\[
y^{4}-3xy^{2}-x^{4}-x^{3}-x-1=0.
\]
Since
\[
y^{4}-3xy^{2}-x^{4}-x^{3}-x-1=(y-x-1)(y+x+1)(y^{2}+x^{2}-x+1),
\]
we see that
\[
x=y-1.
\]
Namely,
\[
\sqrt[3]{\bigg(\frac{\Pi_{q^{3}}}{\Pi_{q^{9}}}\bigg)^{2}-1}=\sqrt{\frac{\Pi_{q}}{\Pi_{q^{9}}}}-1.
\]
Cubing both sides of the above identity and simplifying we can easily
deduce the result (\ref{eq:c2-4}) for $0<q<1$. By analytic continuation,
we know that (\ref{eq:c2-4}) holds for $|q|<1.$ \qed

\section*{Acknowledgements}

 The first author was supported by the Natural Science Basic Research
Plan in Shaanxi Province of China (No. 2017JQ1001), the Initial Foundation
for Scientific Research of Northwest A\&F University (No. 2452015321)
and the Fundamental Research Funds for the Central Universities (No.
2452017170). The second author was supported by the National Natural
Science Foundation of China (Grant No. 11371184) and the Natural Science
Foundation of Henan Province (Grant No. 162300410086, 2016B259, 172102410069).

\end{document}